\documentclass[12pt,a4paper,leqno]{article}

\usepackage{amssymb,exscale}
\usepackage[centertags]{amsmath}
\usepackage{amsthm}

\numberwithin{equation}{section}

\swapnumbers
\theoremstyle{definition}
\newtheorem{theorem}[equation]{Theorem}
\newtheorem{lemma}[equation]{Lemma}

\renewcommand{\phi}{\varphi}

\renewcommand{\(}{\bigl(}
\renewcommand{\)}{\bigr)\vphantom{)}}

\renewcommand{\Re}{\operatorname{Re}}

\newcommand{\BES}{\operatorname{BES}}

\newcommand{\eps}{\varepsilon}
\newcommand{\ga}{\gamma}

\newcommand{\de}{\delta}

\newcommand{\B}{\mathcal B}
\newcommand{\C}{\mathcal C}
\renewcommand{\P}{\mathcal P}
\newcommand{\la}{\lambda}

\newcommand{\Ex}{\mathbb E}
\renewcommand{\Pr}[1]{\,\mathbb P\,\(\,#1\,\)\,}
\newcommand{\R}{\mathbb R}
\newcommand{\CC}{\mathbb C}

\newcommand{\ti}{\tilde}

\newcommand{\sif}{$\sigma$-field}

\newcommand{\One}{\mathbf1}

\newcommand{\equi}{\;\Longleftrightarrow\;}

\def\emailwww#1#2{\par\qquad {\tt #1}\par\qquad {\tt #2}\medskip}

\begin{document}

\title{From random sets\\
  to continuous tensor products:\\
  answers to three questions of W.~Arveson}
\author{Boris Tsirelson}
\date{}

\maketitle

\begin{abstract}
The set of zeros of a Brownian motion gives rise to a product system
in the sense of William Arveson (that is, a continuous tensor product
system of Hilbert spaces). Replacing the Brownian motion with a Bessel
process we get a continuum of non-isomorphic product systems.
\end{abstract}

\section*{Introduction}
``The term \emph{product system} is a less tortured contraction of the
phrase \emph{continuous tensor product system of Hilbert spaces}''
(Arveson \cite[p.~6]{Ar89}). The theory of product systems, elaborated
by W.~Arveson in connection with \mbox{$ E_0 $-semigroups} and quantum fields
(see \cite{Ar96}, \cite{Ar89} and refs therein) suffers from lack of
rich sources of examples. I propose such a source by combining
A.~Vershik's idea of a \emph{measure type factorization}
\cite[Sect.~1c]{TV}, my own idea of a \emph{spectral type of a noise}
\cite[Sect.~2]{TsU}, and J.~Warren's idea (private communication,
Nov.~1999) of constructing a measure type factorization from a given
random set. The new rich source of examples leads to rather simple
answers to three questions of Arveson; see Sections 2,4,5 for the
questions, and Theorems \ref{2.1}, \ref{4.2} and \ref{5.4} for the
answers.

It is interesting to compare measure type factorizations with
so-called \emph{noises} (a less tortured substitute for such phrases as
\emph{homogeneous continuous tensor product system of probability
spaces} or \emph{stationary probability measure factorization}), see
\cite{TV}, \cite{WaPo}, \cite{TsF} and refs therein. Theory of noises
is able to answer two out of the three questions of Arveson, however,
the new approach makes it easier. I still do not know whether the
third question (see Sect.~4) also has a noise-theoretic answer, or
not.

\section{The construction}
Consider the standard Brownian motion $ B(\cdot) $ in $ \R $, and the
random set
\[
Z_{t,a} = \{ s \in [0,t] : B(s) = a \} \, ,
\]
where $ a,t \in (0,\infty) $ are parameters.\footnote{%
 When writing $ Z_{t,a} $ I always assume that $ a,t \in (0,\infty) $
 unless otherwise stated; the reservation applies when I write, say, $
 Z_{\infty,0} $.}
The set $ Z_{t,a} $ may
be treated as a random variable taking on values in the space $ \C_t $
of all closed subsets of $ [0,t] $.\footnote{%
 Also the empty set $ \emptyset $ belongs to $ \C_t $.}
There is a natural Borel \sif\ $ \B_t $ on $ \C_t $, and $ (\C_t,\B_t)
$ is a standard Borel space. Moreover, $ \C_t $ is a compact metric
space w.r.t.\ the Hausdorff metric $ \rho_t ( C_1, C_2 ) = \inf \{
\eps > 0 : C_1 \subset (C_2)_{+\eps} \,\&\, C_2 \subset (C_1)_{+\eps}
\} $ (here $ C_{+\eps} $ means the $ \eps $-neighborhood of $ C $),
and $ \B_t $ is the Borel \sif\ of the metric space $ (\C_t,\rho_t)
$. Let $ P_{t,a} $ be the law of the $ \C_t $-valued random variable $
Z_{t,a} $, then $ (\C_t,\B_t,P_{t,a}) $ is a probability space.

\begin{lemma}\label{1.1}
$ P_{t,a_1} \sim P_{t,a_2} $; that is, measures $ P_{t,a_1} $ and $
P_{t,a_2} $ are equivalent ($=\,$mutually absolutely continuous) for
all $ a_1,a_2 \in (0,\infty) $.
\end{lemma}

\begin{proof}
Consider the random time $ T_a = \min \{ t \in [0,\infty) : B(t) = a
\} $. The shifted set $ Z_{\infty,a} - T_a $ is independent of $ T_a $
and distributed like $ Z_{\infty,0} $. Thus, $ P_{\infty,a} $ is a mix
of shifted copies of $ P_{\infty,0} $, weighted according to the law
of $ T_a $. However, laws of $ T_{a_1},
T_{a_2} $ are equivalent measures, therefore $ P_{\infty,a_1} \sim
P_{\infty,a_2} $, which implies $ P_{t,a_1} \sim P_{t,a_2} $.
\end{proof}

Denote by $ \P_t $ the set of all probability measures on $
(\C_t,\B_t) $ that are equivalent to $ P_{t,a} $ for some (therefore,
all) $ a \in (0,\infty) $. The triple $ (\C_t,\B_t,\P_t) $ is an
example of a structure called \emph{measure-type space.}

Denote by $ P_{s,a} \otimes P_{t,a} $ the law of the random set $ C_1
\cup (C_2+s) $, where $ C_1 \in \C_s $ is distributed $ P_{s,a} $, and
$ C_2 \in \C_t $ is distributed $ P_{t,a} $, and $ C_1, C_2 $ are
independent; of course, $ C_2+s \subset [s,t] $ is the shifted $ C_2
$.

\begin{lemma}\label{1.2}
$ P_{s,a} \otimes P_{t,a} \sim P_{s+t,a} $ for all $ s,t,a \in
(0,\infty) $.
\end{lemma}

\begin{proof}
The conditional distribution of the set $ ( Z_{s+t,a} \cap [s,s+t] ) -
s $, given the set $ Z_{s,a} $, is the mix (over $ x $) of its
conditional distributions, given $ Z_{s,a} $ and $ B(s) = x $. The
latter conditional distribution, being equal to $ P_{t,|a-x|} $, belongs
to $ \P_t $ (except for $ x=a $, which case may be
neglected). Therefore the former conditional distribution also belongs
to $ \P_t $.
\end{proof}

We cannot identify the Cartesian product $ \C_s \times \C_t $ with $
\C_{s+t} $, since natural maps $ \C_{s+t} \to \C_s \times \C_t $ and $
\C_s \times \C_t \to \C_{s+t} $ are not mutually inverse (in fact,
both are non-invertible). However, $ \P_{s+t} \{ C : s \in C \} = 0
$;\footnote{%
 I mean, of course, that $ P \( \{ C \in \C_{s+t} : s \in C \} \) = 0
 $ for some (therefore all) $ P \in \P_{s+t} $.}
neglecting some sets of probability $ 0 $, we get
\begin{equation}\label{1a}
(\C_s,\B_s,\P_s) \otimes (\C_t,\B_t,\P_t) =
(\C_{s+t},\B_{s+t},\P_{s+t}) \, ,
\end{equation}
or simply $ \P_s \otimes \P_t = \P_{s+t} $ for $ s,t \in (0,\infty)
$.

In order to introduce Hilbert spaces $ L_2 (\C_t,\B_t,\P_t) $ note
that Hilbert spaces $ L_2 (\C_t,\B_t,P_1) $ and $ L_2 (\C_t,\B_t,P_2)
$ for $ P_1,P_2 \in \P_t $ are in a natural unitary correspondence;
namely, $ \psi_1 \in L_2 (\C_t,\B_t,P_1) $ corresponds to $ \psi_2 \in
L_2 (\C_t,\B_t,P_2) $ if
\[
\psi_2 = \sqrt{ \frac{P_1}{P_2} } \psi_1 \, ,
\]
where $ \frac{P_1}{P_2} $ is the Radon-Nikodym density. Define an
element $ \psi $ of $ L_2 (\C_t,\B_t,\P_t) $ as a family $ \psi =
(\psi_P)_{P\in\P_t} $ satisfying $ \psi_P \in L_2 (\C_t,\B_t,P) $ and
\[
\psi_{P_2} = \sqrt{ \frac{P_1}{P_2} } \psi_{P_1} \quad \text{for all }
P_1,P_2 \in \P_t \, .
\]
Clearly, $ L_2 (\C_t,\B_t,\P_t) $ is a separable Hilbert space,
naturally isomorphic to every $ L_2 (\C_t,\B_t,P) $, $ P \in \P_t
$.\footnote{%
 Intuitively we may think that $ \sqrt P \psi_P = \psi $ for all $ P
 \in \P_t $. See also \cite{Ac}.}
Relation \eqref{1a} gives
\begin{equation}\label{1b}
L_2 (\C_s,\B_s,\P_s) \otimes L_2 (\C_t,\B_t,\P_t) = L_2
(\C_{s+t},\B_{s+t},\P_{s+t})
\end{equation}
in the sense that the two Hilbert spaces are \emph{naturally}
isomorphic.

However, \eqref{1b} is only a part of requirements stipulated in the
definition of a product system \cite[Def.~1.4]{Ar89}. The point is
that \eqref{1b} holds for each $ (s,t) $ individually; nothing was
said till now about measurability in $ s,t $. In order to get a
product system, we need a \emph{measurable} unitary correspondence
between spaces $ L_2 (\C_t,\B_t,\P_t) $ for different $ t $, making
the map implied by \eqref{1b} jointly measurable. The correspondence
need not be natural, but our case is especially nice, having a natural
correspondence described below.

For every $ \la \in (0,\infty) $ the random process $ t \mapsto
\sqrt\la B(t/\la) $ is a Brownian motion, again. Therefore the two
random sets $ \{ s : B(s) = a \} $ and $ \{ s : \sqrt\la B(s/\la) = a
\} = \la \cdot \{ s : B(s) = a / \sqrt\la \} $ are identically
distributed. It means that the ``rescaling'' map $ R_\la : \C_1 \to
\C_\la $, defined by $ R_\la (C) = \la \cdot C $, sends $
P_{1,a/\sqrt\la} $ to $ P_{\la,a} $. Accordingly, it sends $ \P_1 $ to
$ \P_\la $. We define a unitary operator $ \ti R_t : L_2
(\C_1,\B_1,\P_1) \to L_2 (\C_t,\B_t,\P_t) $ by
\[
( \ti R_t \psi )_{R_t(P)} (R_t(C)) = \psi_P (C) \quad \text{for $ P
  $-almost all $ C \in \C_1 $},
\]
for all $ \psi \in L_2 (\C_1,\B_1,\P_1) $ and $ P \in \P_1 $; of
course, $ R_t (P) $ is the $ R_t $-image of $ P $ (denoted also by $ P
\circ R_t^{-1} $). The disjoint union $ E = \cup_{t\in(0,\infty)} L_2
(\C_t,\B_t,\P_t) $ (not a Hilbert space, of course) is now
parametrized by the Cartesian product $
(0,\infty) \times L_2 (\C_1,\B_1,\P_1) $, namely, $ (t,\psi) \in
(0,\infty) \times L_2 (\C_1,\B_1,\P_1) $ parametrizes $ \ti R_t (\psi) \in
L_2 (\C_t,\B_t,\P_t) \subset E $. We equip $ E $ with the Borel
structure that corresponds to the natural Borel structure on $
(0,\infty) \times L_2 (\C_1,\B_1,\P_1) $. Linear operations and the
scalar product are Borel measurable (on their domains) for trivial
reasons. It remains to consider the multiplication $ E \times E \to E
$,
\[
E \times E \supset H_s \times H_t \ni (\psi_1,\psi_2) \mapsto \psi_1
\otimes \psi_2 \in H_s \otimes H_t = H_{s+t} \subset E \, ;
\]
it must be Borel measurable.\footnote{%
 I do not distinguish between $ H_s \otimes H_t $ and $ H_{s+t} $ in
 the notation. A cautious reader may insert a notation for the natural
 unitary operator $ H_s \otimes H_t \to H_{s+t} $.}
In other words, we consider $ \psi = \ti R^{-1}_{s+t} \( \ti R_s(\psi_1)
\otimes \ti R_t(\psi_2) \) $ as an $ H_1 $-valued function of four
arguments $ s,t \in (0,\infty) $, $ \psi_1,\psi_2 \in H_1 $; we have
 to check that the function is jointly Borel measurable. After
 substituting all relevant definition it boils down to $ C =
 R^{-1}_{s+t} \( (R_s C_1) \cup (s+R_t C_2) \) $ treated as a $ \C_1
 $-valued  function of four arguments $ s,t \in (0,\infty) $, $
 C_1,C_2 \in \C_1  $; the reader may check that the function is
 jointly Borel measurable. So, Hilbert spaces
\[
H_t = L_2 ( \C_t, \B_t, \P_t )
\]
form a product system.

\section{Units}
Every measure $ P \in \P_t $ has an atom, since $ \Pr{ Z_{t,a} =
\emptyset } > 0 $; in fact, $ \{ \emptyset \} $ is the only atom of $
P $.

For every $ t \in (0,\infty) $ the space $ H_t = L_2 (\C_t,\B_t,\P_t)
$ contains a special element $ v_t $ defined by
\[
(v_t)_P (C) = \begin{cases}
  \frac1{ \sqrt{ P(\{\emptyset\}) } } & \text{if $ C = \emptyset $},\\
  0 & \text{otherwise}.
\end{cases}
\]
Clearly, $ v_{s+t} = v_s \otimes v_t $ for all $ s,t \in (0,\infty)
$. Also, $ \| v_t \| = 1 $ for all $ t $.

A unit of a product system $ (H_t) $ is a family $
(u_t)_{t\in(0,\infty)} $ such that $ u_t \in H_t $ for all $ t \in
(0,\infty) $, and $ u_s \otimes u_t = u_{s+t} $ for all $ s,t \in
(0,\infty) $, and the map $ \R \ni t \mapsto u_t \in \cup_t H_t $ is
measurable, and $ u_t \ne 0 $ for some $ t $ (which implies $ u_t \ne
0 $ for all $ t $); see \cite[p.~10]{Ar96}, \cite[Sect.~4]{Ar89}.

The family $ (v_t) $ is a unit, since $ \ti R_t^{-1} (v_t) $ is
measurable in $ t $; in fact, it is constant, $ \ti R_t^{-1} (v_t) =
v_1 $.

If $ (u_t) $ is a unit (of a product system) then $ (e^{i\la t} u_t) $
is also a unit for every $ \la \in \CC $. All these units may be
called equivalent. Some product systems contain non-equivalent
units. Some product systems contain no units at all. The trivial
product system (consisting of one-dimensional Hilbert spaces) contains
a unit, and all its units are equivalent. Arveson \cite[p.~12]{Ar96}
asked: is there a nontrivial product system that contains a unit but
does not contain non-equivalent units? The product system constructed
in Sect.~1 appears to be such an example; the question is answered by
the following result. (Note however that the question is already
answered by noise theory; I mean the system of \cite[Sect.~5]{TV}.)

\begin{theorem}\label{2.1}
Every unit $ (u_t) $ is of the form $ u_t = e^{i\la t} v_t $.
\end{theorem}

\begin{proof}
Every $ \psi \in H_t $ determines a measure $ |\psi|^2 $ on $
(\C_t,\B_t) $ by\footnote{%
 Do not confuse the \emph{measure} $ |\psi|^2 $ with the \emph{number}
 $ \| \psi \|^2 $, the squared norm; in fact, $ \| \psi \|^2 =
 (|\psi|^2) (\C_t) $, the total mass.}
\begin{equation}\label{2.2}
\frac{ |\psi|^2 }{ P } = |\psi_P|^2 \quad \text{for some (therefore,
  all) } P \in \P_t \, .
\end{equation}
Note that $ | \psi_1 \otimes \psi_2 |^2 = |\psi_1|^2 \otimes
|\psi_2|^2 $ whenever $ \psi_1 \in H_s $, $ \psi_2 \in H_t $. If $
(u_t) $ is a unit, then $ |u_s|^2 \otimes |u_t|^2 = |u_{s+t}|^2 $. We
may assume that $ \| u_t \| = 1 $ for all $ t $ (since $ (u_t/\|u_t\|)
$ is a unit equivalent to $ (u_t) $, see \cite[Th.~4.1]{Ar89}), then $
|u_t|^2 $ is a probability measure. Applying \cite[Th.~4.1]{Ar89}
again we get $ \langle u_t, v_t \rangle = e^{\ga t} $ for some $ \ga
\in \CC $. However, for every $ \psi \in H_t $
\begin{gather*}
\langle \psi, v_t \rangle = \int \psi_P \overline{ (v_t)_P } \, dP =
  \psi_P (\emptyset) \frac1{ \sqrt{ P(\{\emptyset\}) } }
  P(\{\emptyset\}) \, , \\
| \langle \psi, v_t \rangle |^2 = | \psi_P (\emptyset) |^2
  P(\{\emptyset\}) = |\psi|^2 (\{\emptyset\}) \, .
\end{gather*}
Applying it to $ \psi = u_t $ we get $ |u_t|^2 (\{\emptyset\}) =
e^{2\Re\ga t} $. In combination with the property $ |u_s|^2 \otimes
|u_t|^2 = |u_{s+t}|^2 $ it shows that $ |u_t|^2 $ is the law of the
Poisson point process with intensity $ (-2\Re\ga) $ on $ [0,t]
$.\footnote{%
 A simple way to check it: divide $ (0,t) $ into $ n $ equal
 intervals; each of them is free of $ C $ (distributed $ |u_t|^2 $)
 with probability $ e^{2\Re\ga t/n} $, independently of
 others. Consider $ n = 2, 4, 8, 16, \dots $}
Thus, $ |u_t|^2 $ is concentrated on finite sets $ C \in \C_t $. On
the other hand, being absolutely continuous w.r.t.\ $ \P_t $, the
measure $ |u_t|^2 $ is concentrated on sets $ C \in \C_t $ with no
isolated points. Therefore $ |u_t|^2 $ is concentrated on $ C =
\emptyset $ only. It means that $ \Re\ga = 0 $, that is, $ \ga = i \la
$, $ \la \in \R $. So, $ \| u_t \| = 1 $, $ \| v_t \| = 1 $ and $
\langle u_t, v_t \rangle = e^{i\la t} $; therefore $ u_t = e^{i\la t}
v_t $.
\end{proof}

\section{Using Bessel processes}
Introduce a parameter $ \de \in (0,2) $ and consider the random set
\[
Z_{t,a,\de} = \{ s \in [0,t] : \BES_{\de,a} (s) = 0 \} \, ,
\]
and its law $ P_{t,a,\de} $;
here $ \BES_{\de,a} (\cdot) $ is the Bessel process of dimension $
\de $ started at $ a $ (see \cite[Chap.~XI, Defs 1.1 and
1.9]{RY}). As before, $ t,a \in (0,\infty) $. The law $ P_{t,a,1} $
of $ Z_{t,a,1} $ is equal to the law $ P_{t,a} $ of $ Z_{t,a} $ of
Sect.~1, since $ \BES_{1,a} $ is distributed like $ |B(\cdot)+a|
$. The structure of $ Z_{\infty,0,\de} $ was well-understood long
ago;\footnote{%
 Namely, $ Z_{\infty,0,\de} $ is the closure of the range of a stable
 subordinator of index $ 1-\de/2 $ (see \cite[Example 6]{PY}); it is
 of Hausdorff dimension $ 1-\de/2 $ near every point \cite{BG}.}
especially, measures $ P_{t,0,\de_1} $ and $ P_{t,0,\de_2} $ for $
\de_1 \ne \de_2 $ are mutually singular. Measures $ P_{t,a,\de_1} $
and $ P_{t,a,\de_2} $ (where $ a > 0 $) are not singular because of a
common atom ($ Z_{t,a,\de} = \emptyset $ with a positive probability).

Below, $ \mu \ll \nu $ means that a measure $ \mu $ is absolutely
continuous w.r.t.\ a measure $ \nu $; $ \mu \sim \nu $ means $ \mu
\ll \nu \,\&\, \nu \ll \mu $.

\begin{lemma}\label{3.1}
(a) $ P_{t,a_1,\de} \sim P_{t,a_2,\de} $;

(b) if $ \de_1 \ne \de_2 $, $ \mu \ll P_{t,a,\de_1} $ and $ \mu \ll
P_{t,a,\de_2} $, then $ \mu $ is concentrated on $ \{ \emptyset \} $.
\end{lemma}

\begin{proof}
Similarly to the proof of Lemma \ref{1.1}, consider the random time $
T_a = \min \{ s \in [0,\infty) : \BES_{\de,a} (s) = 0 \} $; $ T_a \in
(0,\infty) $ almost sure (since $ \de < 2 $). The shifted set $
Z_{\infty,a,\de} - T_a $
is independent of $ T_a $ and distributed like $ Z_{\infty,0,\de}
$. Statement (a) follows from the fact that laws of $ T_{a_1}, T_{a_2}
$ are equivalent measures. Statement (b): $ \mu $ is concentrated on
sets that must have two different Hausdorff dimensions near each
point; the only such set is $ \emptyset $.
\end{proof}

\begin{lemma}
$ P_{s,a,\de} \otimes P_{t,a,\de} \sim P_{s+t,a,\de} $.
\end{lemma}

The proof is quite similar to the proof of Lemma \ref{1.2}.

The Bessel process has the same scaling property as the Brownian
motion: the process $ t \mapsto \sqrt\la \BES_{\de,a/\sqrt\la} (t/\la)
$ has the law $ P_{t,a,\de} $ irrespective of $ \la \in (0,\infty) $.

So, all properties of Brownian motion, used in Sect.~1, hold for
Bessel processes. Generalizing the construction of Sect.~1 we get a
product system $ (H_{t,\de})_{t\in(0,\infty)} $ for every $ \de \in
(0,2) $. The product system of Sect.~1 corresponds to $ \de = 1 $.

\section{Continuum of non-isomorphic product systems}
\hfill\parbox{10cm}{%
``At this point, we are not even certain of the \emph{cardinality} of $
\Sigma \, $! It is expected that $ \Sigma $ is uncountable, but this
has not been proved.''\hfill W.~Arveson \cite[p.~12]{Ar96}.
}

An isomorphism between two product systems $ (H_t) $, $ (H'_t) $ is
defined naturally as a family $ (\theta_t)_{t\in(0,\infty)} $ of
unitary operators $ \theta_t : H_t \to H'_t $ such that, first, $
\theta_{s+t} (\psi_1 \otimes \psi_2) = \theta_s (\psi_1) \otimes
\theta_t (\psi_2) $ whenever $ \psi_1 \in H_s $, $ \psi_2 \in H_t $,
and second, $ \theta_t (\psi) $ is jointly measurable in $ t $ and $
\psi $; see \cite[p.~6]{Ar89}. Are there uncountably many
non-isomorphic product systems? This question, asked by Arveson
\cite[p.~12]{Ar96}, will be answered here in the positive by showing
that product systems $ (H_{t,\de}) $ for different $ \de $ are
non-isomorphic.

Consider the projection operator (the index $ \de $ is suppressed)
\[
Q_t : H_t \to H_t \, , \qquad (Q_t \psi)_P (C) =
 \begin{cases} \psi_P (C) & \text{if $ C = \emptyset $}, \\
   0 & \text{otherwise},
 \end{cases}
\]
just the orthogonal projection onto the one-dimensional subspace
corresponding to the atom of $ \P_{t,\de} $. Given $ 0 < r < s < t $,
we introduce an operator $ Q_{t,(r,s)} = Q_r \otimes \One_{s-r}
\otimes Q_{t-s} $ on the space $ H_t = H_r \otimes H_{s-r} \otimes
H_{t-s} $; of course, $ \One_{s-r} $ is the identical operator on $
H_{s-r} $. Operators $ Q_{t,E} $ are defined similarly for every
elementary set (that is, a union of finitely many intervals) $ E
\subset (0,t) $.\footnote{%
 For example, $ Q_{t,(r,s)\cup(u,v)} = Q_r \otimes \One_{s-r} \otimes
 Q_{u-s} \otimes \One_{v-u} \otimes Q_{t-v} $ for $ 0 < r < s < u < v
 < t $; $ Q_{t,(0,s)} = \One_s \otimes Q_{t-s} $; $ Q_{t,(s,t)} = Q_s
 \otimes \One_{t-s} $; $ Q_{t,(0,t)} = \One_t $; $ Q_{t,\emptyset} =
 Q_t $.}
Clearly,
\[
(Q_{t,E} \psi)_P (C) = 
 \begin{cases} \psi_P (C) & \text{if $ C \subset E $}, \\
   0 & \text{otherwise}.
 \end{cases}
\]
Note a relation to measures $ |\psi|^2 $ defined by \eqref{2.2}:
\begin{equation}\label{4.1}
\langle Q_{t,E} \psi, \psi \rangle = |\psi|^2 ( \{ C \in \C_t : C
\subset E \} ) \, .
\end{equation}

\begin{theorem}\label{4.2}
If $ \de_1 \ne \de_2 $ then product systems $ (H_{t,\de_1}) $, $
(H_{t,\de_2}) $ are non-isomorphic.
\end{theorem}

\begin{proof}
Assume the contrary: operators $ \theta_t : H_{t,\de_1} \to
H_{t,\de_2} $ are an isomorphism of the product systems. The system $
(H_{t,\de_1}) $ has a unit, and all its units are equivalent, which is
Theorem \ref{2.1} when $ \de_1 = 1 $, and a (straightforward)
generalization of Theorem \ref{2.1} for arbitrary $ \de_1 $. The same
for the other product system $ (H_{t,\de_2}) $. It follows that
operators $ Q_t $ are preserved by isomorphisms; $ Q_t \theta_t =
\theta_t Q_t $ (that is, $ Q_t^{(\de_2)} \theta_t = \theta_t
Q_t^{(\de_1)} $). Tensor products of these operators are also
preserved:
\[
Q_{t,E} \theta_t = \theta_t Q_{t,E} \, .
\]
In combination with \ref{4.1} it gives for $ \psi \in H_{t,\de_1} $
\begin{equation}\label{4.3}
|\psi|^2 (A) = | \theta_t \psi |^2 (A)
\end{equation}
for every $ A $ of the form $ A = A_E = \{ C \in \C_t : C \subset E \}
$ where $ E $ is an elementary set. However, $ A_{E_1\cap E_2} =
A_{E_1} \cap A_{E_2} $, and the \sif\ generated by sets $ A_E $ is the
whole $ \B_t $. It follows (by Dynkin Class Theorem) that \ref{4.3}
holds for all $ A \in \B_t $, that is,
\[
| \psi |^2 = | \theta_t \psi |^2 \quad \text{for all } \psi \in
H_{t,\de_1} \, ,
\]
which contradicts to Lemma \ref{3.1}(b).
\end{proof}

\section{Asymmetry via countable random sets}
The law $ P_{t,a} $ of the random set $ Z_{t,a} $ of Sect.~1 is
asymmetric in the sense that $ P_{t,a} $ is not invariant under the
time reversal
\[
\C_t \ni C \mapsto t - C \in \C_t
\]
(of course, $ t-C = \{ t-s : s \in C \} $). However, the measure type
$ \P_t $ is symmetric; therefore the product system $ (H_t) $ is
symmetric, which means existence of unitary operators $ \theta_t : H_t
\to H_t $ such that, first, $ \theta_{s+t} ( \psi_1 \otimes \psi_2 ) =
\theta_t (\psi_2) \otimes \theta_s (\psi_1) $ whenever $ \psi_1 \in
H_s $, $ \psi_2 \in H_t $, and second, $ \theta_t (\psi) $ is jointly
measurable in $ t
$ and $ \psi $; see \cite[p.~12]{Ar96}, \cite[p.~6]{Ar89}. It was
noted by Arveson \cite[p.~6]{Ar89} that we do not know if an
arbitrary product system is symmetric. Apparently, the first
example of an asymmetric product system is ``the noise made by a
Poisson snake'' of J.~Warren \cite{WaPo}; there, asymmetry emerges
from a random countable closed set that has points of accumulation
from the left, but never from the right. A different, probably simpler
way from such sets to asymmetric product systems is presented here.

Our first step toward a suitable countable random set is choosing a
(nonrandom) set $ S \subset [0,\infty) $ and a function $ \la : S
\times S \to [0,\infty) $ such that

(a) $ S $ is closed, countable, $1$-periodic (that is, $ s \in S
\equi s+1 \in S $ for $ s \in [0,\infty) $), totally ordered (that
is, no strictly decreasing infinite sequences), $ 0 \in S $, and $ S
\cap (0,1) $ is infinite;\footnote{\label{exS}
 An example: $ S = \{ k - 2^{-l} : k,l = 1,2,3,\dots \} \cup \{
 0,1,2,\dots \} $; another example: $ S = \{ k - 2^{-l} - 2^{-l-m} :
 k,l,m = 1,2,3,\dots \} \cup \{ k - 2^{-l} : k,l = 1,2,3,\dots \} \cup
 \{ 0,1,2,\dots \} $.}

(b) $ \la ( s_1, s_2 ) > 0 $ whenever $ s_1, s_2 \in S $, $ s_1 < s_2
\le s_1+1 $; and $ \la ( s_1, s_2 ) = 0 $ whenever $ s_1, s_2 \in S $ do
not satisfy $ s_1 < s_2 \le s_1 + 1 $;

(c) denoting by $ s_+ $ the least element of $ S \cap (s,\infty) $ we
have
\[
\la ( s, s_+ ) = \frac1{s_+-s} \, , \qquad \sum_{s'\in S, \, s'>s_+}
\la(s,s') \le 1
\]
for all $ s \in S $.

On the second step we construct a Markov process $ \( X(t)
\)_{t\in[0,\infty)} $ that jumps, from one point of $ S $ to another,
according to the rate function $ \la(\cdot,\cdot) $. Initially, $ X(0)
= 0 $. We introduce independent random variables $ \tau_{0,s} $ for $
s \in S \cap (0,1] $ such that $ \Pr{ \tau_s > t } = e^{-\la(0,s)t} $
for all $ t \in [0,\infty) $. We have $ \inf_s \tau_s > 0 $, since $
\sum_s \la (0,s) < \infty $. We let
\[
X(t) = 0 \text{ for } t \in [0,T_1) \, , \qquad
X(T_1) = s_1 \, ,
\]
where random variables $ T_1 \in (0,\infty) $ and $ s_1 \in S $ are
defined by
\[
T_1 = \inf_s \tau_s = \tau_{s_1} \, .
\]
The first transition of $ X(\cdot) $ is constructed. Now we construct
the second transition, $ X(T_2-) = s_1 $, $ X(T_2) = s_2 $ using rates
$ \la(s_1,s) $; and so on. It may happen (in fact, it happens almost
always) that $ \sup_k T_k = T_\infty < \infty $, and then (almost
always) $ X(T_k) \to s_\infty \in S $ (recall that $ S $ is closed).
We let $ X(T_\infty) = s_\infty $ and construct the next transition
of $ X(\cdot) $ using rates $ \la(s_\infty,s) $. And so on, by a
transfinite recursion over countable ordinals, until exhausting the
time domain $ [0,\infty ) $. Almost surely, $ X(t) \in S $ is well-defined
for all $ t \in [0,\infty) $, and $ X(t) \to \infty $ for $ t \to
\infty $.

The last step is simple. We define the random set $ Z_{\infty,0,S} $
as the closure of the set of all instants when $ X(\cdot) $
jumps. That is, $ Z_{\infty,0,S} $ is the set of all $ t $ such that $
X(t-\eps) < X(t+\eps) $ for all $ \eps \in (0,t) $. Instead of
starting at $ 0 $ we may start at another point $ a \in S $, which
leads to another process $ X_a (\cdot) $ and random set $
Z_{\infty,a,S} $; the law $ P_{t,a,S} $ of $ Z_{t,a,S} =
Z_{\infty,a,S} \cap [0,t] $ is a probability measure on $ (\C_t,\B_t)
$.

\begin{lemma}
$ P_{t,a_1,S} \sim P_{t,a_2,S} $ for all $ a_1, a_2 \in S $.
\end{lemma}

\begin{proof}
(Similar to \ref{1.1}.)
Consider the random time $ T_a = \min Z_{a,S} $, just the instant of
the first jump: $ X_a (T_a-) = a $, $ X_a (T_a) > a $. The conditional
distribution of the shifted set (without the first point), $
(Z_{\infty,a,S}-T_a) \setminus \{0\} $, given $ T_a $ and $ X_a (T_a)
$, is $ P_{\infty,X_a(T_a),S} $. Thus, $ P_{\infty,a,S} $ is a mix of
shifted copies of $ P_{\infty,b,S} \cup \{ 0 \} $ for various $ b \in
S \cap (a,a+1] $. However, $ P_{\infty,b,S} = P_{\infty,b+1,S} $ for 
all $ b \in S $. It remains to note that the joint law of $ T_{a_1} $
and $ \( X_{a_1} (T_{a_1}) \bmod 1 \) $ is equivalent to the joint law
of $ T_{a_2} $ and $ \( X_{a_2} (T_{a_2}) \bmod 1 \) $.
\end{proof}

We denote by $ \P_{t,S} $ the set of all probability measures on $
(\C_t,\B_t) $ that are equivalent to $ P_{t,a,S} $ for some
(therefore, all) $ a \in S $.

\begin{lemma}\label{5.2}
$ P_{s,a,S} \otimes P_{t,a,S} \sim P_{s+t,a,S} $ for all $ s,t \in
(0,\infty) $, $ a \in S $.
\end{lemma}

\begin{proof}
(Similar to \ref{1.2}.)
The conditional distribution of the set $ ( Z_{s+t,a,S} \cap [s,s+t] )
- s $, given the set $ Z_{s,a,S} $, is the mix (over $ b $) of its
conditional distributions, given $ Z_{s,a,S} $ and $ X_a(s) = b $. The
latter conditional distribution, being equal to $ P_{t,b,S} $, belongs
to $ \P_{t,S} $. Therefore the former conditional distribution also
belongs to $ \P_{t,S} $.
\end{proof}

Now we can construct the corresponding product system $
(H_{t,S})_{t\in[0,\infty)} $ as before. Though, scaling invariance is
absent; unlike Sect.~1, $ R_t $ does not send $ \P_{1,S} $ to $
\P_{t,S} $. We have no \emph{natural} correspondence between spaces $
L_2 (\C_t,\B_t,\P_{t,S}) $, but still, \emph{some} Borel-measurable
correspondence exists; I do not dwell on this technical issue.

A more important point: in contrast to previous sections, the product
system $ (H_{t,S}) $ contains non-equivalent units (since the law of a
Poisson point process on $ (0,t) $ is absolutely continuous w.r.t.\ $
\P_{t,S} $). Unlike Sect.~4, an isomorphism need not preserve
projection operators $ Q_t $ and measures $ |\psi|^2 $, which prevents
us from deriving asymmetry of the product system $ (H_{t,S}) $ just
from asymmetry of measure types $ \P_{t,S} $. Instead, we'll adapt
some constructions of \cite{TsF} (see (2.15) and (3.4) there).

As before, $ Q_t : H_{t,S} \to H_{t,S} $ is the one-dimensional
projection operator corresponding to the atom $ \{ \emptyset \} $ of $
\P_{t,S} $ (you see, $ \Pr{ Z_{t,a,S} = \emptyset } > 0 $). Introduce
operators
\[
U_{t,p,n} = \( (1-p) Q_{t/n} + p \One_{t/n} \)^{\otimes n}
\]
on $ H_t = H_{t/n} \otimes \dots \otimes H_{t/n} = H_{t/n}^{\otimes n}
$ (here $ p \in (0,1) $ is a parameter).\footnote{%
 Of course, $ \One_t $ is the identical operator on $ H_{t,S} $.}
It is just multiplication by a
function of $ C \in \C_t $; the function counts intervals $ (\frac k
n, \frac{k+1}n ) $ that contain points of $ C $, and returns $ p^m $
where $ m $ is the number of such intervals. For $ n \to \infty $,
operators $ U_{t,p,n} $ converge (in the strong operator topology) to
\[
U_{t,p} = \lim_{n\to\infty} U_{t,p,n} \, , \qquad (U_{t,p}\psi)_P (C)
= p^{|C|} \psi_P (C) \, ,
\]
just multiplication by $ p^{|C|} $ where $ |C| $ is the cardinality of
$ C $; naturally, $ p^{|C|} = 0 $ for infinite sets $ C $. (In fact, $
U_{t,p_1} U_{t,p_2} = U_{t,p_1p_2} $.) The operator $ U_{t,1-} =
\lim_{p\to1-} U_{t,p} $ is especially interesting:
\[
(U_{t,1-}\psi)_P (C) = \begin{cases}
  \psi_P (C) & \text{if $ C $ is finite}, \\
  0 & \text{otherwise}.
\end{cases}
\]
(In fact, $ U_{t,1-} $ is the projection onto the stable ($ = \, $
linearizable) part of the product system \cite[(2.15)]{TsF}, which is
not used here.)

Operators $ U_{t,p} $ correspond to a particular unit (or rather,
equivalence class of units) of the product system $ (H_{t,S})
$. However, we may do the same for any given unit $ u = (u_t)
$. Namely,
\begin{gather*}
Q_{t,u} \psi = \frac{ \langle \psi, u_t \rangle }{ \langle u_t, u_t
  \rangle } u_t \quad \text{for } \psi \in H_t \, ; \\
U_{t,p,n,u} = \( (1-p) Q_{t/n,u} + p \One_{t/n} \)^{\otimes n} \, ; \\
U_{t,p,u} = \lim_{n\to\infty} U_{t,p,n,u} \, .
\end{gather*}
Existence of the limit is an easy matter, since operators $
U_{t,p,n,u} $ for all $ n $ belong to a single commutative
subalgebra. Even simpler, we may take $ \lim_{n\to\infty}
U_{t,p,2^n,u} $, the limit of a \emph{decreasing} sequence of
commuting operators.

\begin{lemma}\label{5.3}
$ U_{t,1-,u} = U_{t,1-} $ for all units $ u $ of the product system $
(H_{t,S}) $.
\end{lemma}

\begin{proof}
Let $ u = (u_t) $ and $ v = (v_t) $ be two units; we'll prove that $
U_{t,1-,u} = U_{t,1-,v} $. Due to \cite[Th.~4.1]{Ar89} we may assume
that $ \| u_t \| = 1 $, $ \| v_t \| = 1 $ and $ \langle u_t, v_t
\rangle = e^{-\ga t} $ for some $ \ga \in [0,\infty) $. An elementary
calculation (on the plane spanned by $ u_t, v_t $) gives\footnote{%
 It is not about product systems, just two vectors in a Hilbert
 space.}
\[
\| Q_{t,u} - Q_{t,v} \| = \sqrt{ 1 - e^{-2\ga t} } \, .
\]
Opening brackets in $ U_{t,p,n,u} = \( (1-p) Q_{t/n,u} + p \One_{t/n}
\)^{\otimes n} $ we get a sum of $ 2^n $ terms, each term being a
tensor product of $ n $ factors. After rearranging the factors (which
changes the term, of course, but does not change its norm), a term
becomes simply $ (1-p)^k p^{n-k} Q_{\frac k n t,u} \otimes
\One_{\frac{n-k}n t} $. We see that
\[
\| U_{t,p,n,u} - U_{t,p,n,v} \| \le \Ex \| Q_{\frac k n t,u} -
Q_{\frac k n t,v} \| \, ,
\]
where the expectation is taken w.r.t.\ a random variable $ k $ having
the binomial distribution $ \operatorname{Bin} (n,1-p) $. Using
concavity of $ \sqrt{1-e^{-2\ga t}} $ in $ t $,
\[
\Ex \| Q_{\frac k n t,u} - Q_{\frac k n t,v} \| = \Ex \sqrt{ 1 -
  e^{-2\ga kt/n} } \le \sqrt{ 1 - e^{-2\ga \Ex kt/n}} = \sqrt{ 1 -
  e^{-2\ga t(1-p) }} \, ,
\]
therefore
\begin{gather*}
\| U_{t,p,n,u} - U_{t,p,n,v} \| \le \sqrt{ 1 - e^{-2\ga t(1-p) }} \quad
  \text{for all } n \, ; \\
\| U_{t,p,u} - U_{t,p,v} \| \le \sqrt{ 1 - e^{-2\ga t(1-p) }} \, ;
\end{gather*}
so, $ \| U_{t,1-,u} - U_{t,1-,v} \| = 0 $.
\end{proof}

Informally, the distinction between empty and non-empty sets $ C \in
\C_t $ is relative (to a special unit) and non-invariant (under
isomorphisms of product systems), while the distinction between finite
and infinite sets $ C \in \C_t $ is absolute, invariant.

For any $ C \in \C_t $ denote by $ C' $ the set of all accumulation
points of $ C $; clearly, $ C' \in \C_t $, and $ C' = \emptyset $ if
and only if $ C $ is finite. We proceed similarly to Sect.~4, but $ C'
$ is used here instead of $ C $. Given an elementary set $ E \subset
(0,t) $, we define operators $ Q'_{t,E} $ by
\[
\( Q'_{t,E} \psi \)_P (C) = \begin{cases}
  \psi_P (C) & \text{if $ C' \subset E $}, \\
  0 & \text{otherwise}.
\end{cases}
\]
We do not worry about boundary points of $ E $, since $ \P_{t,S}
$-almost all $ C $ avoid them. Operators $ Q'_{t,E} $ are tensor
products of operators $ U_{s,1-} $. (For example, if $ E = (r,s)
$, $ 0 < r < s < t $, then $ Q'_{t,E} = U_{r,1-} \otimes \One_{s-r}
\otimes U_{t-s,1-} $.) By Lemma \ref{5.3}, every isomorphism
preserves $ U_{s,1-} $; therefore it preserves $ Q'_{t,E} $. Given $
\psi \in H_{t,S} $, we define a measure $ {|\psi|'}^2 $ on $
(\C_t,\B_t) $ as the image of the measure $ |\psi|^2 $ (defined by
\eqref{2.2}) under the map $ \C_t \ni C \mapsto C' \in \C_t
$. Similarly to \eqref{4.3} we see that $ {|\psi|'}^2 $ is preserved
by isomorphisms (even though $ |\psi|^2 $ is not).

\begin{theorem}\label{5.4}
If $ S'' \ne \emptyset $ then the product system $ (H_{t,S}) $ is
asymmetric.\footnote{%
 Of course, $ S'' $ means $ (S')' $; recall examples of $ S $ on page
 \pageref{exS}.}
\end{theorem}

\begin{proof}
Assume the contrary: the product system is symmetric; $ \theta_t :
H_{t,S} \to H_{t,S} $, $ \theta_{s+t} ( \psi_1 \otimes \psi_2 ) =
\theta_t (\psi_2) \otimes \theta_s (\psi_1) $ for $ \psi_1 \in H_{s,S}
$, $ \psi_2 \in H_{t,S} $. Then
\[
\theta_t Q'_{t,E} = Q'_{t,t-E} \theta_t \, .
\]
It follows that
\begin{equation}\label{5.5}
R_t \( {|\psi|'}^2 \) = {| \theta_t \psi |'}^2 \quad \text{for } \psi
\in H_{t,S} \, ;
\end{equation}
here $ R_t ( {|\psi|'}^2 ) $ is the image of the measure $ {|\psi|'}^2 $
under the time reversal $ R_t : \C_t \to \C_t $, $ R_t (C) = t - C
$. However, for $ \P_{t,S} $-almost all $ C \in \C_t $, $ C $ is
totally ordered, therefore $ C' $ is also totally ordered. Both
measures, $ {|\psi|'}^2 $ and $ {|\theta_t \psi|'}^2 $, being
absolutely continuous w.r.t.\ $ \P_{t,S} $, are concentrated on
totally ordered sets. In combination with \eqref{5.5} it means that
they are concentrated on finite sets. So, $ C'' = \emptyset $ for $
\P_{t,S} $-almost all $ C \in \C_t $.

The Markov process $ X(\cdot) $ consists of ``small jumps'' $ X(t) =
\( X(t-) \)_+ $ and ``big jumps'' $ X(t) > \( X(t-) \)_+ $.\footnote{%
 As before, $ s_+ $ is the least element of $ S \cap (s,\infty) $.}
The rate of big jumps never exceeds $ 1 $. The rate of small jumps
results in the mean speed $ 1 $ in the sense that $ X(t) - t $ is a
martingale between big jumps. There is a chance that $ X(\cdot) $
increases by $ 1 $ (or more) by small jumps only (between big
jumps). In such a case, $ S'' \ne \emptyset $ implies $ Z''_{t,a,S} \ne
\emptyset $. So, $ \{ C \in \C_t : C'' \ne \emptyset \} $ is not $
\P_{t,S} $-negligible, in contradiction to the previous paragraph.
\end{proof}

\bigskip
\filbreak
\begingroup
{
\small
\begin{sc}
\parindent=0pt\baselineskip=12pt

School of Mathematics, Tel Aviv Univ., Tel Aviv
69978, Israel
\emailwww{tsirel@math.tau.ac.il}
{http://www.math.tau.ac.il/$\sim$tsirel/}
\end{sc}
}
\filbreak

\endgroup

\end{document}